\documentclass[12pt]{article}

\usepackage{amsfonts,amssymb,amsbsy,amsmath,amscd}

\newtheorem{theorem}{Theorem}[section]
\newtheorem{proposition}{Proposition}[section]
\newtheorem{lemma}{Lemma}[section]
\newtheorem{tcorollary}{Corollary}[theorem]
\newtheorem{pcorollary}{Corollary}[proposition]

\def\proc#1{\medbreak\noindent{\it #1}\hspace{1ex}\ignorespaces}
\def\ep{\noindent{\hfill $\Box$}}

\begin{document}

\title{Sturmian words, $\beta$-shifts, and transcendence
}

\author{Dong Pyo\ Chi  and  DoYong\ Kwon
    \thanks{The first author was supported by
    Korea Research Foundation Grant (KRF-2000-015-DP0031),
    and the second author was supported by BK21 Program in Korea
    and by KIAS Research Fund (No.03-0155-001).}}
\maketitle

\makeatletter
\renewcommand\section
   {\@startsection{section}{1}{20pt}{\baselineskip}{-1pc}{\bfseries}}
\makeatother

\makeatletter
\renewcommand{\@seccntformat}[1]{\csname the#1\endcsname. \hspace{0.05em}}
\makeatother
\begin{abstract}
Consider the minimal $\beta$-shift containing the shift space
generated by given Sturmian word. In this paper we characterize
such $\beta$ and investigate its combinatorial, dynamical and
topological properties and prove that such $\beta$ are
transcendental numbers.
\end{abstract}
%
%
%
%       I  N  T  R  O  D  U  C  T  I  O  N
%
%
%
\section{Introduction.}
Sturmian words are infinite words over binary alphabet $A$
whose factors of length $n$ are exactly $n+1$. It is known that
Sturmian words are aperiodic infinite words with minimal
complexity \cite{CH}. These critical words admit some equivalent
definitions in different manners. They can be coded from the
irrational billiards on a unit square or equivalently from the
irrational rotation on $\mathbb{R}/\mathbb{Z}$ under a certain
partition \cite{MH}. And they also carry the balanced properties
which will be defined later \cite{MH}.\\
\indent Let $\beta>1$ be a real number. We consider
$\beta$-transformation $T_\beta$ on $[0,1]$ defined by $T_\beta :
x\mapsto\beta x\mod 1$. Then the $\beta$-expansion of $x\in[0,1]$,
denoted by $d_\beta (x)$, is a sequence of integers determined by
the following rule: $$d_\beta (x)=(x_i)_{i\geq1}\ \mathrm{if\ and\
only\ if\ }x_i = \lfloor \beta T_\beta^{i-1} (x) \rfloor,$$ where
$\lfloor t \rfloor$ is the largest integer not greater than $t$.
If $\beta$ is not an integer, $d_\beta (x)$ is a sequence over the
alphabet $A=\{ 0,1, \cdots , \lfloor \beta \rfloor \}$. When
$\beta$ is an integer, the digits $x_i$ belong to $A=\{ 0,1,
\cdots , \beta-1\}$. The $\beta$-shift $S_\beta$ is the closure of
$\{ d_\beta (x)|x\in[0,1) \}$ with respect to the topology of
$A^\mathbb{N}$. In \cite{Pa}, Parry completely characterized
$S_\beta$ in terms of $d_\beta (1)$ and the lexicographic order
given in $A^\mathbb{N}$. From Parry's result we note that the
collection of $\beta$-shifts is totally ordered. The main concern
of this article is about the minimal $\beta$-shift containing the
shift space generated by a
Sturmian word.\\
\indent For $1<\beta<2$, we call $\beta$ a Sturmian number if the
set $\{T_\beta ^n 1\}_{n\geq0}$ is infinite and $$1-\frac1 \beta
\leq T_\beta ^n 1 \leq 1\ \ \mathrm{for\ all\ } n\geq0.$$ We show
that if $S_\beta$ minimally contains the shift space generated by
some Sturmian word then $\beta$ is a Sturmian number. And
conversely we also show that $d_\beta (1)$ is a Sturmian word for
every Sturmian number $\beta$. The definition of Sturmian number
is naturally generalized for $\beta\geq2$. Then the alphabet
involved is extended from $\{0,1\}$ to $\{a,b\}$ with $0\leq
a<b=\lfloor \beta\rfloor$. This gives a large class of specified
$\beta$-transformations $T_\beta$ and moreover for such $\beta$,
the diameters of closure of $\{T_\beta ^n 1\}_{n\geq0}$ are
minimal in a certain sense. We also prove the transcendence of
such $\beta$. That is a partial answer to the question posed by
Blanchard \cite{Bl}.

%
%
%
%       P R E L I M I N A R I E S
%
%
%
\section{Preliminaries.}
In this section we briefly review the terminology on words.\\
\indent An \textit{alphabet} is a finite set $A$ which may have an
order or not. We consider the free monoid $A\sp *$ generated by
$A$.  We mean by a \textit{word}, any element in $A\sp *$, whereas
a \textit{(right) infinite word} is one in $A\sp \mathbb{N}$.
$A\sp *$ has a natural binary operation called the
\textit{concatenation} and the \textit{empty word} $\varepsilon$
serves as its identity element. For a word $w$, we denote by $|w|$
the \textit{length} of $w$, i.e., if $w=a_1 a_2 \cdots a_n$ and
all $a_i \in A$, then $|w|=n$.
Sometimes a word of length $1$ is called a \textit{letter}.\\
\indent A (finite) word $w$ is said to be a \textit{factor} (a
\textit{prefix}, a \textit{suffix}, respectively) of a finite or
infinite word $u$ provided $u$ is expressed as $u=xwy$ ($u=wy$,
$u=xw$, respectively) for some words $x$ and $y$ with $x$ finite.
The set of factors of a word $x$ is denoted by $F(x)$. Among
factors of $x$, we mean by $F_n (x)$, the set of factors with
length $n$. For $X\subset A\sp *$ or $X\subset A\sp \mathbb{N}$,
we define
$$F(X) := \bigcup_{x \in X} F(x).$$ $F_n (X)$ is defined in a
similar way. We say a subset $X$ of $A\sp *$ is \textit{factorial}
if all factors of its elements
are again in $X$, that is, $x\in X$ implies $F(x) \subset X.$\\
\indent The \textit{reversal} of a word $w=a_1 a_2 \cdots a_n$,
where $a_1, a_2, \cdots ,a_n$ are letters, is the word
$\tilde{w}=a_n a_{n-1} \cdots a_1$. And a \textit{palindrome word}
is a word $w$ such that $w=\tilde{w}$. In
particular the empty word and all letters are palindrome words.\\
\indent The \textit{complexity function} of an infinite word $x$
is, for each integer $n\geq0$, the cardinality of the set of
factors of length $n$ in $x$. In formula, one can write
$$P(x,n):= \mathrm{Card} (F_n (x)).$$ For a subset $X$ of $A\sp
\mathbb{N}$, $P(X,n)$ can be defined similarly. The
\textit{frequency} of a word is a partial function that measures
the relative occurrence of $w$ in an infinite word $x$. Suppose
$\mu_x^N (w)$ is the number of occurrences of $w$ in the prefix of
length $N+|w|-1$ of $x$. In other words, if $w=a_1 a_2 \cdots a_m$
and $x=x_1 x_2 \cdots$, then we have
$$\mu_x^N (w)=\mathrm{Card}\left(\{j|x_{1+j}=a_1, \cdots,
x_{m+j}=a_m, 0\leq j\leq N-1\}\right).$$ Then the frequency of $w$
in $x$ is $$\mu_x (w)= \lim_{N \to \infty}\frac1 N \mu_x^N (w),\ \
\mathrm{if\ the\ limit\
exists}.$$\\
\indent If an alphabet has an order, it can be extended to $A\sp
*$ and $A\sp \mathbb{N}$ lexicographically which is called the
\textit{lexicographic order}. More precisely, suppose $A$ is an
ordered alphabet. Given $x, y \in A\sp \mathbb{N}$, we denote
$x<y$ if there exist a nonnegative integer $m$ such that $x_1=y_1,
\cdots, x_{m-1}= y_{m-1}$ and $x_{m} < y_{m}$. Here $x_i, y_j$'s
are all letters. For finite words $x=x_1 \cdots x_m$ and $y=y_1
\cdots y_n$, we write $x<y$ if $x_1 \cdots x_m 0 0 \cdots < y_1
\cdots y_n 0 0
\cdots$.\\
\indent The set $A\sp \mathbb{N}$ is well endowed with a metric in
a sense that the metric generates the usual product topology of
$A\sp \mathbb{N}$. For any $x, y \in A\sp \mathbb{N}$, we define
the distance between $x$ and $y$ by $d(x,y)=2^{-n}$, where $n=\min
\{k\geq 0| x_k \neq y_k \}$.

%
%
%     S T U R M I A N   W O R D S    A N D
%
%     L E X I C O G R A P H I C   O R D E R
%
%
\section{Sturmian words and lexicographic order.}

A \textit{Sturmian word} is an infinite word $s$ with its
complexity function $P(s,n)=n+1$ for any nonnegative integer $n$.
Since $P(s,1)=2$, Sturmian words are forced to be infinite words
over the alphabet $\{0,1\}$ by renaming if necessary. Thus we
assume $A=\{0,1\}$ unless stated explicitly.

\proc{Example 1.} The Fibonacci word is an infinite word defined
by
$$f_0 =0,\ f_1 =01,\ f_{n+2} =f_{n+1} f_n,\ n\geq 0.$$
We note the sequence of $|f_n |$ is the famous integer sequence of
Fibonacci numbers. By recursive arguments, we get
$$f= \lim_{n \to \infty} f_n = 0100101001001010010100100101001001
\cdots.$$ One can also find that $f$ is Sturmian. See Chapter 2 of
Lothaire's book \cite{Lo}. \medbreak

\indent In the present section, we give two alternative
characterizations of Sturmian words. The first one is the
`balanced property.'\\
\indent The \textit{height} $h(x)$ of a word $x$ is the number of
the occurrences of $1$ in $x$. We say a subset $X$ of $A\sp *$ is
\textit{balanced} if for any $x, y \in X,$ $|h(x)-h(y)|\leq 1$
whenever $|x|=|y|.$ An infinite word $s$ is also called
\textit{balanced} if $F(s)$ is balanced.  Morse and Hedlund
\cite{MH} showed the following.

\begin{theorem} \label{balanced}
Suppose $s$ is an infinite word. Then $s$ is Sturmian if and only
if $s$ is aperiodic and balanced.
\end{theorem}

In \cite{CH}, Coven and Hedlund described the balanced property in
more detail.

\begin{theorem} \label{palindrome}
Let $X$ be a factorial subset of $A\sp *$. Then $X$ is unbalanced
if and only if there exists a palindrome word $w$ such that both
$0w0$ and $1w1$ lie in $X$.
\end{theorem}

\indent The \textit{slope} of a nonempty word $x$ is the real
number $\pi(x)=h(x)/|x|$. Let $x$ be an infinite balanced word,
and $x_n$ be the prefix of $x$ with length $n\geq 1$. Then we see
that the sequence $(\pi(x_n))_{n\geq 1}$ converges to some value
$\pi(x)$(see \cite{Lo}), which we will call the \textit{slope} of
$x$.

\proc{Example 2.} The slope of Fibonacci word can be easily
computed. Let $F_n =|f_n|$. Then $h(f_n)=F_{n-2}$ and
$$\pi(f)=\lim_{n \to \infty} \frac{F_{n-2}}{F_n}
=\frac{1}{\tau^2},$$ where $\tau= (1+\sqrt{5})/2$.
\medbreak

\indent Next we survey Sturmian words in a different point of
view.\\
\indent For a real number $t$, $\lceil t\rceil$ the smallest
integer not less than $t$, and $\{t\}$ is the fraction part of
$t$, i.e., $t=\lfloor t \rfloor + \{t\}$. Let $\alpha$, $\rho$ be
two real numbers with $0\leq \alpha\leq 1$. We now define two
infinite words over $\{0, 1\}$. Consider, for nonnegative integer
$n$,
$$s_{\alpha,\rho}(n) = \lfloor \alpha(n+1) + \rho \rfloor -
\lfloor \alpha n +\rho \rfloor,$$
$$s'_{\alpha,\rho}(n) = \lceil \alpha(n+1) + \rho \rceil -
\lceil \alpha n +\rho \rceil.$$ The infinite words
$s_{\alpha,\rho}$, $s'_{\alpha,\rho}$ are termed a \textit{lower
mechanical word} and an \textit{upper mechanical word}
respectively with \textit{slope} $\alpha$ and \textit{intercept}
$\rho$. And \textit{mechanical words} refer to either lower or
upper mechanical words. Until now we have introduced the
terminology `slope' twice for infinite words. But two definitions
of slope coincide for infinite balanced words.
 One notes that if $\rho$ and $\rho'$ differ by an integer, then
$s_{\alpha,\rho}=s_{\alpha,\rho'}$ and
$s'_{\alpha,\rho}=s'_{\alpha,\rho'}$ hold. Hence with no loss of
generality, we assume $0\leq \rho <1$ through the article unless
stated explicitly. If $\alpha n +\rho$ is not an integer for any
$n\geq 0$, we have $s_{\alpha,\rho}=s'_{\alpha,\rho}$. Otherwise
in the case that $\alpha n +\rho$ is an integer for some $n>0$, we
get
$$s_{\alpha,\rho}(n-1)=1,\quad s_{\alpha,\rho}(n)=0,$$
$$s'_{\alpha,\rho}(n-1)=0,\quad s'_{\alpha,\rho}(n)=1.$$
Thus if $\alpha$ is irrational, then $s_{\alpha,\rho}$,
$s'_{\alpha,\rho}$ are the same possibly except only one factor of
length at most $2$. Worthy of comment is the special case where $\rho$ has the value $0$.
If $\alpha$ is irrational, we see
$$s_{\alpha,0}=0c_\alpha, \quad s'_{\alpha,0}=1c_\alpha$$ for some
infinite word $c_\alpha$. Here the word $c_\alpha$ is called the
\textit{characteristic word} of $\alpha$. Morse and Hedlund
\cite{MH} also characterized Sturmian words in terms of mechanical
words.

\begin{theorem}
Suppose $s$ is an infinite word. Then $s$ is Sturmian if and only
if $s$ is irrationally mechanical.
\end{theorem}

\indent The following propositions are another known results on
Sturmian words that we will use later.

\begin{proposition}[\cite{Mi}] \label{factor}
For two Sturmian words $s$, $t$, if they have the same slope, then
$F(s)=F(t)$. And $F(s)\cap F(t)$ is finite otherwise.
\end{proposition}

\begin{proposition}[\cite{He}]
The shift space generated by a Sturmian word is minimal.
\end{proposition}

\begin{proposition}[\cite{Pv}] \label{morphism}
A map $E$ on $\{0,1\}$ is defined by $E(0)=1$, $E(1)=0$. And
consider the natural extension to $\{0,1\}\sp \mathbb{N}$ by
concatenations. Then $E(c_\alpha)=c_{1-\alpha}$.
\end{proposition}

\indent If $s$ is some Sturmian word, we denote by
$\overline{\mathcal{O}}(s)$ the shift space generated by $s$,
i.e., the orbit closure of $s$. For the shift map we write
$\sigma$. One can note that
$$\sigma^n (s_{\alpha,\rho})= s_{\alpha,\{\alpha n +\rho\}}.$$ The
same is true for upper mechanical words. It is well known that
irrational rotations on a circle are ergodic and their orbits are
all dense.

\begin{proposition}\label{orbit}
Let $s$ be a Sturmian word with slope $\alpha$. Then
$\overline{\mathcal{O}}(s)$ is the set of all mechanical words of
slope $\alpha$.
\end{proposition}

The proof is a consequence of a lemma.

\begin{lemma}
For fixed irrational $\alpha \in (0,1)$, $s_{\alpha,\rho}$ is
continuous from the right and $s'_{\alpha,\rho}$ from the left as
functions of $\rho$.
\end{lemma}
\proc{Proof.} Let $\epsilon>0$, $s_{\alpha,\rho_0}$,
$s'_{\alpha,\rho_0}$ be given. We choose an integer $N>0$ such
that $2^{-N}<\epsilon$. And put
$$\delta_1= \min \{ 1-\{\alpha n +\rho_0 \}\ |\ 0\leq n \leq N+1 \}.$$
Then $0 \leq \rho-\rho_0 < \delta_1 /2$ implies
$d(s_{\alpha,\rho},s_{\alpha,\rho_0})<\epsilon.$ For the upper
mechanical word, we define $\delta_2$ by the minimum of nonzero
fractions $\{\alpha n +\rho_0 \}$ for $0\leq n \leq N+1$. Then it
is true that $0 \leq \rho_0 -\rho < \delta_2 /2$ implies
$d(s'_{\alpha,\rho},s'_{\alpha,\rho_0})<\epsilon.$ If $\rho_0 =0$,
then we can assume $\rho_0 =1$. \ep \medbreak

\proc{Proof of Proposition \ref{orbit}.} By the minimality of
$\overline{\mathcal{O}}(s)$ we may assume $s=s'_{\alpha,0} =
1c_\alpha$. Thus one has $\sigma^n (s)=s'_{\alpha,\{\alpha n\}}$.
Since $\alpha$ is irrational, $\alpha n$ is never an integer for
nonzero $n$. Hence $$\sigma^n (s)=s_{\alpha,\{\alpha
n\}}=s'_{\alpha,\{\alpha n\}}$$ holds for any $n\geq 1$.\\
\indent Given $s_{\alpha,\rho}$ and $s'_{\alpha,\rho}$ we can pick
two increasing sequences of integers $(p_n )_{n\geq0}$, $(q_n
)_{n\geq0}$ such that $\{\alpha p_n \} \searrow \rho$ and
$\{\alpha q_n \} \nearrow \rho$. This is possible because the set
of $(\{\alpha n \})_{n\geq0}$ is dense in $[0,1]$. Then one finds
$$ \lim_{n \to \infty} \sigma^{p_n} (s)=s_{\alpha,\rho},
\ \lim_{n \to \infty} \sigma^{q_n} (s)=s'_{\alpha,\rho}.$$ \indent
Conversely assume $t \in \overline{\mathcal{O}}(s)$. We have to
show that $t$ is a Sturmian word of slope $\alpha$. In fact, $t$
is balanced since $F(t)\subset F(s)$. And the minimality implies
the aperiodicity of $t$.
Thus $t$ is a Sturmian word.\\
\indent $t$ has the slope $\alpha$. Indeed if it has a different
slope, then $F(s) \cap F(t)$ is finite by Proposition
\ref{factor}. But every factor of $t$ also occurs in $s$.
 \ep \medbreak

\indent In the consecutive section, what we need critically is the
lexicographic orders between Sturmian words. At first we know

\begin{proposition} \label{lexicographic}
Suppose $\alpha \in (0,1)$ is irrational and $\rho, \rho'
\in[0,1)$ are real. Then $$s_{\alpha,\rho}<s_{\alpha,\rho'}\
\textit{if\ and\ only\ if\ } \rho<\rho'.$$
\end{proposition}
\proc{Proof.} See \cite{Lo}. \ep \medbreak

\indent From this proposition, we derive

\begin{theorem}
Let $\alpha$ be an irrational number in $(0,1)$. Then $$0c_\alpha
< s_{\alpha,\rho} <1c_\alpha\ \textit{for\ any\ } 0<\rho<1.$$
\end{theorem}

\proc{Proof.} If $s_{\alpha,\rho} (0)=0$, then
$\{\alpha+\rho\}=\alpha+\rho$. Since $c_\alpha =
s_{\alpha,\alpha}$ and $\sigma
(s_{\alpha,\rho})=s_{\alpha,\alpha+\rho}$, it follows from
Proposition \ref{lexicographic} that $0c_\alpha <s_{\alpha,\rho}$.
For the case $s_{\alpha,\rho} (0)=1$, we get
$\{\alpha+\rho\}=\alpha+\rho-1$. Thus $\sigma (s_{\alpha,\rho}) =
s_{\alpha,\alpha+\rho-1}$. Noting $\alpha+\rho-1<\alpha$, we see
$s_{\alpha,\rho}<1c_\alpha$.
\ep \medbreak

\indent The next corollary indeed was known to Borel and Laubie
\cite{BL}.

\begin{tcorollary} \label{admissible}
Let $s$ be a Sturmian word with a slope $\alpha$. Then $1c_\alpha$
is the maximal element and $0c_\alpha$ the minimal one in
$\overline{\mathcal{O}}(s)$. In particular we have $$1c_\alpha
>\sigma^n (1c_\alpha) \ \ \textit{and}\ \ 0c_\alpha<\sigma^n
(0c_\alpha)$$ for all positive integers $n$.
\end{tcorollary}

\indent Given a Sturmian word $s$, the theorem above presents an
explicit algorithm to find the maximal and minimal points of
$\overline{\mathcal{O}}(s)$.

\proc{Example 3.} Let $f$ be the Fibonacci word. For any
$x\in\overline{\mathcal{O}}(f)$, we know
$$0c_{\tau^{-2}} \leq x\leq 1c_{\tau^{-2}}\ \ \mathrm{or}$$
$$00100101001001010010100100\cdots \leq x \leq
  10100101001001010010100100\cdots .$$
And the inequality is the best possible.
\medbreak

\indent If two Sturmian words have different slopes, the
lexicographic order is given by the following theorem.

\begin{theorem} \label{order of characteristics}
For two irrational $\alpha$, $\beta$ in $(0,1)$, we have
$$c_\alpha <c_\beta \ \ \textit{if\ and\ only\ if\ }\ \alpha<\beta.$$
\end{theorem}

\proc{Proof.} For a lower mechanical word $s=s_{\alpha,\rho}$, the
height of a prefix $u=s(0)s(1)\cdots s(n-1)$ is given by the
integer $$h(u)=\lfloor \alpha n+ \rho\rfloor.$$ This fact
guarantees the proof.\\
\indent Suppose $c_\alpha <c_\beta$ or equivalently $s_{\alpha,0}
<s_{\beta,0}$. Then one has $\lfloor \alpha n \rfloor < \lfloor
\beta n \rfloor$ for some integer $n\geq 1$, which implies $\alpha
<\beta$. For $\alpha <\beta$, we choose the smallest positive
integer $n$ such that $\lfloor \alpha n \rfloor < \lfloor \beta n
\rfloor$ holds. We find then that $$s_{\alpha,0}=0a_1 \cdots
a_{n-2} 0\cdots < 0a_1 \cdots a_{n-2} 1\cdots =s_{\beta,0}.$$  \ep
\medbreak

%
%
%
%     B E T A - S H I F T S
%     A N D  S T U R M I A N  N U M B E R S
%
%
%
\section{$\beta$-shifts and Sturmian numbers.}

Recall a $\beta$-shift $S_\beta$ is the closure of all
$\beta$-expansions of real numbers in $[0,1)$. Just as the number
1 dominates any number in $[0,1)$, so does $d_\beta(1)$ in
$S_\beta$ with respect to the lexicographic order. Parry \cite{Pa}
showed:

\begin{theorem}
Given $\beta>1$, let $s$ be an element of
$\{0,1,\cdots,\lfloor\beta\rfloor\}\sp \mathbb{N}$.\\
If $d_\beta(1)$ is not finite, then $s$ belongs to $S_\beta$ if
and only if
$$\sigma^n(s)\leq d_\beta (1)\ \ \textit{for\ all\ } n\geq 0.$$
If $d_\beta(1)=d_1\cdots d_m 00\cdots$, then $s$ belongs to
$S_\beta$ if and only if
$$\sigma^n(s)\leq d_1\cdots d_{m-1}(d_m -1)d_1\cdots d_{m-1}(d_m
-1)d_1\cdots\ \textit{for\ all\ } n\geq0.$$
\end{theorem}

Moreover Parry also characterized sequences that can be
$\beta$-expansions of $1$ for some $\beta>1$. Such sequences obey
the next rule.

\begin{theorem} \label{dbeta=1}
A sequence $s\in \{0,1,\cdots,\lfloor\beta\rfloor\}\sp \mathbb{N}$
is a $\beta$-expansion of 1 for some $\beta$ if and only if
$\sigma^n(s)<s$ for $n\geq 1$. In the case, such $\beta$ is
unique.
\end{theorem}

The following proposition that is also due to Parry furnishes a
total order to the collection of all $\beta$-shifts.

\begin{proposition} \label{order of dbeta=1}
Suppose $\beta$, $\gamma\ >1$. Then
$$\beta<\gamma\ \ \textit{if\ and\ only\ if\ \ }
d_\beta(1)<d_\gamma(1).$$
\end{proposition}

\begin{pcorollary}
If $\beta<\gamma$, then $S_\beta \subset S_\gamma$.
\end{pcorollary}

With the above results in mind, we consider the minimal
$\beta$-shift containing the orbit closure of a Sturmian word.\\
\indent Let $s$ be a Sturmian word of slope $\alpha$. Then any
$t\in \overline{\mathcal{O}}(s)$ lies between $0c_\alpha$ and
$1c_\alpha$. For this we use the notation as $0c_\alpha \leq
\overline{\mathcal{O}}(s) \leq 1c_\alpha$. This notation
represents $\beta$-shift as $0^\infty \leq S_\beta \leq
d_\beta(1)$. In both cases, two inequalities are the best
possible. By Theorem \ref{dbeta=1} and Corollary \ref{admissible},
there exists a unique $\beta\in(1,2)$ such that
$d_\beta(1)=1c_\alpha$. Thus we deduce that such $\beta$-shift is
the minimal one we are searching for. Moreover the closure of
$\{\sigma^n(d_\beta(1) )\}_{n\geq0}$ is equal to
$\overline{\mathcal{O}}(s)$ and $0c_\alpha$, $1c_\alpha$ are
accumulation points by Proposition \ref{orbit}. Here the minimal
point informs us that $d_\beta(1-1/\beta)=0c_\alpha$. We state
these facts as a theorem.

\begin{theorem}
Let $s$ be a Sturmian word of slope $\alpha$. If $S_\beta$ is the
smallest $\beta$-shift containing $\overline{\mathcal{O}}(s)$,
then $\overline{ \{\sigma^n(d_\beta(1) )\}}_{n\geq0}  =
\overline{\mathcal{O}}(s)$ and $\beta$ is the unique solution of
$$ 1= \sum_{n=0}^{\infty} \frac{s'_{\alpha,0}(n)}{x^{n+1}}.$$
\end{theorem}

In terms of $\beta$-transformation $T_\beta$, the theorem means

\begin{tcorollary} \label{Sturmian number}
For $\beta$ appearing in the theorem, the closure of $\{T_\beta^n
1\}_{n\geq0}$ is contained in $[1-1/\beta,1]$. Moreover
$1-1/\beta$ and $1$ are accumulation points.
\end{tcorollary}

If no order is given to $A=\{a,b\}$, then one can determine an
order on $A$ so that the minimal $\beta$-shift has a smaller
topological entropy. For a Sturmian word $s$ with slope $\alpha$,
Proposition \ref{morphism} implies
$$0c_{1-\alpha}\leq E(s)\leq 1c_{1-\alpha}.$$
Hence by Theorem \ref{order of characteristics} and Proposition
\ref{order of dbeta=1}, all we have to know is whether $\alpha$ is
smaller than $1/2$ or not.

\proc{Example 4.} Suppose $f$ is an infinite word defined by
$f_0=a$, $f_1=ab$, $f_{n+2}=f_{n+1}f_{n},\ n\geq0$. Then $a=0$,
$b=1$ makes the minimal $\beta$-shift have a smaller topological
entropy than the case $a=1$, $b=0$ since $\tau^{-2}<1/2$.
\medbreak

Because of Corollary \ref{Sturmian number}, we have a good reason
to name such numbers after Sturm.

\proc{Definition.} $\beta\in(1,2)$ is a \textit{Sturmian number}
if $d_\beta (1)$ is aperiodic and $1-1/\beta\leq T_\beta^n 1
\leq1$ for any $n\geq0$. \medbreak

\indent Corollary \ref{Sturmian number} reads as follows: If
$S_\beta$ minimally contains all Sturmian words of the same slope,
then $\beta$ is a Sturmian number. We will next prove the
converse.

\begin{theorem}
If $\beta$ is a Sturmian number, then $d_\beta(1)$ is the Sturmian
word $1c_\alpha$ for some $\alpha$.
\end{theorem}

\proc{Proof.} Given a Sturmian number $\beta$, $d_\beta(1)$ is
aperiodic. Thus it suffices to show $d_\beta(1)$ is balanced. If
it is unbalanced, then Theorem \ref{palindrome} guarantees the
existence of a palindrome word $w$ such that both $0w0$, $1w1$ are
factors of $d_\beta(1)$. For $d_\beta(1)=1d_1 d_2 \cdots$, one
sees $d_\beta(1-1/\beta)=0d_1 d_2 \cdots$. Since $\beta$ is
Sturmian, we get
$$0d_1 \cdots d_n d_{n+1}\leq 0w0<1w1 \leq 1 d_1 \cdots d_n
d_{n+1},$$ where $n$ is the length of $w$. Thus $w$ must be equal
to $d_1 \cdots d_n$, that is,
$$wd_{n+1}\leq w0<w1\leq wd_{n+1}.$$
But this yields a contradiction. Hence $d_\beta(1)$ is a Sturmian
word of some slope $\alpha$ and it dominates all its shifts, and
therefore $d_\beta(1)=1c_\alpha$. \ep \medbreak

\indent In the proof of the theorem above, a crucial step is from
the combinatorial argument. But there is something to declare
about pure real dynamics.

\begin{tcorollary}
For any $\beta\in[1,2]$, either $\overline{\{T_\beta^n 1\}}_{n
\geq 0}$ is finite or has a diameter not less than $1/\beta$.
\end{tcorollary}

%
%
%     D Y N A M I C S   OF   S T U R M I A N
%
%     B E T A - T R A N S F O R M A T I O N S
%
%
\section{Dynamics of Sturmian $\beta$-transformations.}
From now on we extend the alphabet $A=\{0,1\}$ to $\{a,b\}$ with
$0\leq a<b$.\\
\indent In \cite{Bl}, Blanchard suggested the study of real
numbers according to the ergodic properties of their
$\beta$-shifts and classified $\beta$-shifts into five categories.
The terminology on language theory used in the next definition is
referred to \cite{Bl} or the bibliography therein.

\proc{Definition.} The set of real numbers greater than 1 is
categorized into five classes according to their $\beta$-shifts:
\begin{itemize}
  \item $\beta \in \mathcal{C}_1$ if and only if $S_\beta$ is a shift of
  finite type.
  \item $\beta \in \mathcal{C}_2$ if and only if $S_\beta$ is sofic.
  \item $\beta \in \mathcal{C}_3$ if and only if $S_\beta$ is specified.
  \item $\beta \in \mathcal{C}_4$ if and only if $S_\beta$ is synchronizing.
  \item $\beta \in \mathcal{C}_5$ if and only if $S_\beta$ has none of the above properties.
\end{itemize}
\medbreak

From the definition we see the following inclusions.
$$\emptyset \neq \mathcal{C}_1 \subset \mathcal{C}_2 \subset
\mathcal{C}_3 \subset \mathcal{C}_4 \subset (1,\infty),\
\mathcal{C}_5 = (1,\infty)\setminus \mathcal{C}_4 .$$

For any $\beta$ contained in some classes, the morphology of its
$\beta$-expansion $d_\beta (1)$ is totally understood by Parry and
Bertrand-Mathis.

\begin{proposition}[\cite{Pa,BM}]
The following equivalences hold.
\begin{itemize}
  \item $\beta \in \mathcal{C}_1$ if and only if $d_\beta (1)$ is finite.
  \item $\beta \in \mathcal{C}_2$ if and only if $d_\beta (1)$ is ultimately
  periodic.
  \item $\beta \in \mathcal{C}_3$ if and only if there exists $n\in \mathbb{N}$
   such that the number of consecutive $0$'s in $d_\beta (1)$ is less than
   $n$, or equivalently the origin is not an accumulation point of
   $\{T_\beta ^n 1\}_{n\geq0}$.
  \item $\beta \in \mathcal{C}_4$ if and only if some word of $F(S_\beta )$
  does not appear in $d_\beta (1)$, or equivalently $\{T_\beta ^n
  1\}_{n\geq0}$ is not dense in $[0,1]$.
  \item $\beta \in \mathcal{C}_5$ if and only if all words of $F(S_\beta )$
  appear at least once in $d_\beta (1)$, or equivalently $\{T_\beta ^n
  1\}_{n\geq0}$ is dense in $[0,1]$.
\end{itemize}
\end{proposition}

On the other hand Schmeling \cite{Sc} determined each size of the
classes, which is one of the questions asked by Blanchard
\cite{Bl}.

\begin{proposition} \label{size}
$\mathcal{C}_3$ has Hausdorff dimension $1$ and $\mathcal{C}_5$
has full Lebesgue measure.
\end{proposition}

Now we concentrate on a special class of real numbers that is
contained in $\mathcal{C}_3$. For an infinite word $x$,
$\mathrm{alph}(x)$ is the set of letters involved in $x$.

\proc{Definition.} Let $\beta >1$. We call $\beta$ a
\textit{Sturmian number} if $d_\beta (1)$ is a Sturmian word over
a binary alphabet $A=\{a,b\}$, $0\leq a<b=\lfloor\beta\rfloor$. In
particular, $\beta$ is \textit{maximally Sturmian} if it is
Sturmian and $\mathrm{alph} \left(d_\beta (1) \right) =\{ \lfloor
\beta \rfloor -1, \lfloor \beta \rfloor \}$. \medbreak

\proc{Remark.} Not every Sturmian word, of course can be $d_\beta
(1)$. In such cases, $d_\beta (1)=D(1c_\alpha)$ for some
irrational $\alpha \in(0,1)$, where $D$ is a morphism defined by
$D(0)=a$ ($a=\lfloor \beta \rfloor -1$ if $\beta$ is maximally
Sturmian) and $D(1)=b=\lfloor\beta\rfloor$. If $\beta\in(1,2)$,
then $\beta$ is Sturmian in the sense of the previous section,
too. \medbreak

To each $\beta >1$, we associate a real number that is the
diameter of $T_\beta$-orbit of $1$.

\proc{Definition.} $\mathrm{diam}:(1,\infty) \rightarrow [0,1]$ is
a function defined by the diameter of $T_\beta$-orbit of $1$, that
is,
$$\mathrm{diam}(\beta) := \mathrm{diam} \{T_\beta ^n 1\}_{n\geq0}
=\sup \{|x-y|: x,y \in \{T_\beta ^n 1\}_{n\geq0} \}.$$ \medbreak

One can note that if $\beta \in \mathcal{C}_1$ or $\beta \in
(1,\infty)\setminus \mathcal{C}_3$, then $\mathrm{diam}(\beta)=1$
since both $0$ and $1$ lie in the closure of $\{T_\beta ^n
1\}_{n\geq0}$. Thus noteworthy is only the case $\beta \in
\mathcal{C}_3 \setminus \mathcal{C}_1$. We get immediately

\begin{proposition}\label{diameter}
Suppose $\beta$ is Sturmian and $\mathrm{alph}( d_\beta (1)
)=\{a,b\}$ with $0\leq a<b=\lfloor \beta \rfloor$. Then
\begin{itemize}
  \item $\beta \in \mathcal{C}_3 \setminus \mathcal{C}_2$.
  \item $\mathrm{diam}(\beta)= \frac{b-a}{\beta}$.
\end{itemize}
Moreover $\beta$ is maximally Sturmian if and only if $\beta
\notin \mathcal{C}_2$ and $1-1/\beta\leq T_\beta^n 1 \leq1$ for
any $n\geq0$.
\end{proposition}

The above proposition implies that for a maximal Sturmian number
$\beta$, its diameter is minimal in the following sense.

\begin{pcorollary}
For any $\beta >1$, either $\beta \in \mathcal{C}_2$ or
$\mathrm{diam}(\beta)\geq 1/\beta$.
\end{pcorollary}

Proposition \ref{size} shows the set of Sturmian numbers is of
Lebesgue measure zero. Then what about the size of
$\overline{\{T_\beta^n 1\}}_{n \geq 0}$ for a fixed Sturmian
number $\beta$? The last paragraph of this section is devoted to
showing $\overline{\{T_\beta^n 1\}}_{n \geq 0}$ has Lebesgue
measure zero whereas the orbit closure of irrational rotation has
full Lebesgue measure.\\

\indent A $\beta$-transformation $T_\beta$ has an invariant
ergodic measure $\nu_\beta$ whose Radon-Nikodym derivative with
respect to Lebesgue measure is given by
$$ h_\beta(x)= \frac{1}{F(\beta)} \sum_{x<T_\beta^n 1}
\frac{1}{\beta^n},\  \ x\in[0,1].$$ Here $F(\beta)$ is the
normalizing factor that makes $\nu_\beta$ a probability measure.
This measure is known to be the unique measure of maximal entropy
\cite{Ho}. Suppose $d_\beta(1)=\epsilon_0 \epsilon_1 \cdots$.
Parry noted the following formulae:
$$F(\beta)=\int_0^1 \sum_{x<T_\beta^n 1} \frac{1}{\beta^n} dx
=\int_0^1 \left( \sum_{n=0}^\infty
\frac{a_n(x)}{\beta^n}\right)dx,$$ where
$$ a_n(x) = \left\{ \begin{array}{ll}
                  1 & \ \mathrm{if}\  x<T_\beta^n 1, \\
                  0 & \ \mathrm{otherwise,}
                  \end{array} \right. $$
$$F(\beta)=\sum_{n=0}^\infty \frac{1}{\beta^n} \int_0^1 a_n(x) dx
          =\sum_{n=0}^\infty \frac{T_\beta^n 1}{\beta^n}
          =\sum_{n=0}^\infty \frac{(n+1)\epsilon_n}{\beta^{n+1}}.$$
For any $x\in [0,1]$ we are able to compute, if any, the frequency
of $\lfloor\beta\rfloor$ in $d_\beta(x)$. It is $\lim_{n \to
\infty}1/n \sum_{i=0}^{n-1} \chi(T_\beta^i (x))$, where
$$ \chi(x) = \left\{ \begin{array}{ll}
                  0 & \ \mathrm{if}\  0\leq x< \lfloor\beta\rfloor/\beta, \\
                  1 & \ \mathrm{if}\  x\geq \lfloor\beta\rfloor/\beta.
                  \end{array} \right. $$
Owing to the Birkhoff Ergodic Theorem, we can say more. For almost
all $x$ in $[0,1]$, the frequency of $\lfloor\beta\rfloor$ in
$d_\beta(x)$ is equal to
$$\mu_\beta(\lfloor\beta\rfloor)=\frac{1}{F(\beta)} \int_{\frac{\lfloor\beta\rfloor}\beta}^1
\sum_{x<T_\beta^n 1} \frac{1}{\beta^n} dx.$$ The similar reasoning
also applies to the frequency of another digit.

\begin{lemma}
If $\beta$ is Sturmian and $\mathrm{alph}( d_\beta (1) )=\{a,b\}$
with $0\leq a<b=\lfloor \beta \rfloor$, then for almost every $x$
in $[0,1]$, the frequency of $b$ in $d_\beta(x)$ is equal to
\begin{align*}
\mu_\beta(b)
     &=\frac{\mathcal{I}}{F(\beta)}
       =\frac{1}{F(\beta)}\sum_{n=0}^\infty \lceil \alpha n \rceil
        \frac{\epsilon_n}{\beta^{n+1}} \\
     &=\frac{1}{F(\beta)}\left(\sum_{n\in J} \lceil \alpha n \rceil
             \frac{b}{\beta^{n+1}}+\sum_{n\in K} \lceil \alpha n \rceil
             \frac{a}{\beta^{n+1}}\right),
\end{align*}
and the frequency of $a$ in $d_\beta(x)$ is equal to
\begin{align*}
\mu_\beta(a)
     &=\frac{\mathcal{J}}{F(\beta)}
      =\frac{1}{F(\beta)}\left(\sum_{n\in J} \frac{1}{\beta^{n+1}}
      +\sum_{n=0}^\infty
       (n-\lceil\alpha n\rceil) \frac{\epsilon_n}{\beta^{n+1}}\right) \\
     &=\frac{1}{F(\beta)}\left(\sum_{n\in J} \frac{1}{\beta^{n+1}}
       +\sum_{n\in J} (n-\lceil \alpha n \rceil)
        \frac{b}{\beta^{n+1}}
       +\sum_{n\in K} (n-\lceil \alpha n \rceil)
        \frac{a}{\beta^{n+1}}\right),
\end{align*}
where $d_\beta(1)=\epsilon_0 \epsilon_1 \epsilon_2 \cdots$ and
$J=\{n\geq0|\epsilon_n=b\}$, $K=\{n\geq0|\epsilon_n=a\}$.
\end{lemma}

\proc{Proof.} Let $\alpha$ be such number that
$d_\beta(1)=bc_\alpha=\epsilon_0 \epsilon_1 \epsilon_2 \cdots$,
where $c_\alpha$ is the characteristic word of slope $\alpha$, but
$1$ is replaced by $b$ and $0$ by $a$. First we compute
integration involved in the frequency of $b$.
$$\mathcal{I}:=\int_{\frac{\lfloor \beta \rfloor}\beta}^1
               \sum_{x<T_\beta^n 1} \frac{1}{\beta^n} dx
             =\int_{\frac{\lfloor \beta \rfloor}\beta}^1
              \sum_{n=0}^\infty \frac{a_n(x)}{\beta^n} dx
             =\sum_{n=0}^\infty \frac{1}{\beta^n}
             \int_{\frac{\lfloor \beta \rfloor}\beta}^1 a_n(x) dx
             =\sum_{n=0}^\infty \frac{b_n}{\beta^n},$$
where
$$ b_n = \left\{ \begin{array}{ll}
                  T_\beta^n 1-\lfloor \beta \rfloor/\beta & \ \mathrm{if}\  \lfloor \beta \rfloor/\beta\leq T_\beta^n 1, \\
                  0 & \ \mathrm{otherwise.}
                  \end{array} \right. $$
Recall $$T_\beta^n 1=
\frac{\epsilon_n}{\beta}+\frac{\epsilon_{n+1}}{\beta^2}+\cdots,$$
and $\lfloor \beta \rfloor/\beta \leq T_\beta^n 1$ is equivalent
to $\epsilon_n=b$. Then one has
$$\mathcal{I} =\sum_{n\in J} \frac{1}{\beta^n} \left(T_\beta^n 1
               -\frac{\lfloor \beta \rfloor}{\beta}\right)
              =\sum_{n\in J} \sum_{m=n+1}^\infty
              \frac{\epsilon_m}{\beta^{m+1}}.$$
By changing the order of indices, we find
$$\mathcal{I}=\sum_{n=0}^\infty h_n \frac{\epsilon_{n+1}}{\beta^{n+2}}
             =\sum_{n=0}^\infty h_{n-1}
             \frac{\epsilon_n}{\beta^{n+1}},$$
where $h_n$ is the number of $b$ in the word $\epsilon_0
\epsilon_1 \cdots \epsilon_n$ and we put $h_{-1}=0$ by convention.
Noting that $h_n=\lceil \alpha(n+1) \rceil$, we finally get
$$\mathcal{I}=\sum_{n=0}^\infty \lceil \alpha n \rceil
\frac{\epsilon_n}{\beta^{n+1}}.$$

The Birkhoff Ergodic Theorem tells us that for almost every $x$ in
$[0,1]$, the frequency of $a$ in $d_\beta(x)$ is equal to
$$\mu_\beta(a)=\frac{\mathcal{J}}{F(\beta)}=\frac{1}{F(\beta)}
\int_{\frac{a}{\beta}}^{\frac{a+1}{\beta}} \sum_{x<T_\beta^n 1}
\frac{1}{\beta^n} dx.$$ The integration is derived as follows.
\begin{align*}
\mathcal{J}&:=\int_{\frac{a}{\beta}}^{\frac{a+1}{\beta}}
              \sum_{x<T_\beta^n 1} \frac{1}{\beta^n} dx
             =\int_{\frac{a}{\beta}}^{\frac{a+1}{\beta}}
              \sum_{n=0}^\infty \frac{a_n(x)}{\beta^n} dx\\
            &=\sum_{n=0}^\infty \frac{1}{\beta^n}
              \int_{\frac{a}{\beta}}^{\frac{a+1}{\beta}} a_n(x) dx
             =\sum_{n=0}^\infty \frac{b_n}{\beta^n},
\end{align*}
where
$$b_n=
  \begin{cases}
    1/\beta & \text{if $(a+1)/\beta\leq T_\beta^n 1$}, \\
    T_\beta^n 1-a/\beta & \text{if $a/\beta\leq T_\beta^n
    1<(a+1)/\beta$}, \\
    0 & \text{if $T_\beta^n 1< a/\beta$}.
  \end{cases}$$
Since only $a$ and $b$ appear in $d_\beta(1)$, the inequality
$T_\beta^n 1< a/\beta$ never occurs and $(a+1)/\beta\leq T_\beta^n
1$ is reduced to $\lfloor\beta\rfloor/\beta\leq T_\beta^n 1$. Then
the integration is expressed as
$$\mathcal{J} =\sum_{n\in J} \frac{1}{\beta^{n+1}}
               +\sum_{n\in K} \frac{1}{\beta^n} \left(T_\beta^n 1
               -\frac{a}{\beta}\right)
              =\sum_{n\in J} \frac{1}{\beta^{n+1}}+
               \sum_{n\in K} \sum_{m=n+1}^\infty
              \frac{\epsilon_m}{\beta^{m+1}}.$$
By changing the order of indices, we find
$$\mathcal{J}=\sum_{n\in J} \frac{1}{\beta^{n+1}}
              +\sum_{n=0}^\infty (n+1-h_n) \frac{\epsilon_{n+1}}{\beta^{n+2}}
             =\sum_{n\in J} \frac{1}{\beta^{n+1}}
              +\sum_{n=0}^\infty (n-h_{n-1}) \frac{\epsilon_n}{\beta^{n+1}}.$$
Hence we see that
$$\mathcal{J}
     =\sum_{n\in J} \frac{1}{\beta^{n+1}}
      +\sum_{n=0}^\infty
       (n-\lceil\alpha n\rceil) \frac{\epsilon_n}{\beta^{n+1}}.$$
\ep \medbreak

\begin{theorem}
If $\beta$ is Sturmian and $\mathrm{alph}( d_\beta (1) )=\{a,b\}$
with $0\leq a<b=\lfloor \beta \rfloor$, then
$\overline{\{T_\beta^n 1\}}_{n\geq 0}$ is of Lebesgue measure
zero.
\end{theorem}

\proc{Proof.} For any $x$ in $\overline{\{T_\beta^n 1\}}_{n \geq
0}$, the infinite word $d_\beta(x)$ is Sturmian.  The frequency of
$b$ in $d_\beta(x)$, therefore, has the value $\alpha$, while the
frequency of $a$ has $1-\alpha$. We will prove that at least one
of these values is different from those given in the lemma.

At first we suppose $a=0$. Then the integration $\mathcal{I}$ is
given by
$$\mathcal{I}=\sum_{n\in J} \lceil \alpha n \rceil
             \frac{b}{\beta^{n+1}}.$$
For the same reason, one sees
$$\alpha F(\beta)= \sum_{n\in J} \frac{ (\alpha n+\alpha
)b}{\beta^{n+1}}.$$ It holds that $n\in J$ if and only if
$\lceil\alpha n\rceil<\alpha n+\alpha$ because
$\epsilon_n=b(\lceil\alpha(n+1)\rceil-\lceil\alpha n\rceil)$.
Whence
$$\mu_\beta(b)=\frac{\mathcal{I}}{F(\beta)}<\alpha.$$

Next, we suppose $1\leq a<b$ and, in addition, $\lfloor \beta
\rfloor \alpha >1$. We find
$$\alpha F(\beta)=\sum_{n\in J} \frac{(\alpha n+\alpha)b}{\beta^{n+1}}
+\sum_{n\in K} \frac{(\alpha n+\alpha)a}{\beta^{n+1}}.$$ Noting
the index set $J$ contains $0$, let $k_0$ be the smallest element
of $K$. We derive
$$\alpha F(\beta)-\mathcal{I}
=\sum_{n\in J} \frac{(\alpha n+\alpha)-\lceil \alpha n
  \rceil}{\beta^{n+1}}b
 -\sum_{n\in K}
        \frac{\lceil \alpha n \rceil-(\alpha
        n+\alpha)}{\beta^{n+1}}a.$$
We know that $(\alpha n+\alpha)>\lceil \alpha n\rceil$ if $n\in
J$, and $\lceil \alpha n \rceil>(\alpha n+\alpha)$ if $n\in K$.
The series can be written and evolved as
\begin{align*}
\alpha F(\beta)-\mathcal{I}
     &>\sum_{n\in J} \frac{(\alpha n+\alpha)-\lceil \alpha n
        \rceil}{\beta^{n+1}}b
       -\sum_{n\in K}
        \frac{a}{\beta^{n+1}} \\
     &=\left(\frac{\alpha}{\beta} b+\cdots\right)
       -\left(\frac{a}{\beta^{k_0+1}} +\cdots\right) \\
     &>\frac{1}{\beta}-\left(\frac{a}{\beta^{k_0+1}}
        +\cdots\right) > 0,
\end{align*}
because $\lfloor\beta\rfloor \alpha>1$ and $k_0\geq1$. Hence we
have
$$\mu_\beta(b)=\frac{\mathcal{I}}{F(\beta)}<\alpha.$$

If $\lfloor\beta\rfloor \alpha < 1$, then
$$\alpha<\frac{1}{\lfloor\beta\rfloor}\leq \frac{1}{2}.$$
From this we can easily check $k_0=1$. Since $(1-\alpha) F(\beta)$
is represented as
$$(1-\alpha) F(\beta) =\sum_{n\in J} \frac{n+1-(\alpha
n+\alpha)}{\beta^{n+1}}b +\sum_{n\in K} \frac{n+1-(\alpha
n+\alpha)}{\beta^{n+1}}a,$$ we see
\begin{align*}
(1-\alpha) F(\beta)-\mathcal{J}&=\sum_{n\in J}\frac{1+\lceil\alpha
n\rceil-(\alpha n+\alpha)}{\beta^{n+1}}b +\sum_{n\in K}
\frac{1+\lceil \alpha n
 \rceil-(\alpha n+\alpha)}{\beta^{n+1}}a\\
 &\quad \quad -\sum_{n\in J}\frac{1}{\beta^{n+1}}\\
&=\left(\frac{1-\alpha}{\beta} b+\cdots\right)
 +\sum_{n\in K} \frac{1+\lceil \alpha n \rceil-
 (\alpha n+\alpha)}{\beta^{n+1}}a\\
 &\quad \quad -\left(\frac{1}{\beta} +\cdots\right).
\end{align*}
By the assumption, the inequality $(1-\alpha)b>b-1\geq1$ is true.
Hence we have the following inequality:
\begin{align*}
(1-\alpha) F(\beta)-\mathcal{J}
  &>\sum_{n\in K} \frac{1+\lceil \alpha n
 \rceil-(\alpha n+\alpha)}{\beta^{n+1}}a
 -\sum_{n\in J\setminus\{0\}}\frac{1}{\beta^{n+1}}\\
&>\frac{1}{\beta^2} -\sum_{n\in
J\setminus\{0\}}\frac{1}{\beta^{n+1}}>0,
\end{align*}
since the least integer in $J\setminus\{0\}$ is greater than or
equal to $2$ and if $n\in K$, then $1+\lceil\alpha n\rceil-(\alpha
n+\alpha) >1$. We have proved
$$\mu_\beta(a)=\frac{\mathcal{J}}{F(\beta)}<1-\alpha.$$
\ep \medbreak

%
%
%
%     T R A N S C E N D E N C E  O F
%     S T U R M I A N   N U M B E R S
%
%
%
\section{Transcendence of Sturmian numbers.}
We know that $\beta$ is an algebraic integer for all $\beta \in
\mathcal{C}_2$. Then are there transcendental numbers in
$\mathcal{C}_3$, $\mathcal{C}_4$, and $\mathcal{C}_5$? This was
questioned by Blanchard in his paper \cite{Bl}. From Schmeling's
results on the sizes of the classes, $\mathcal{C}_5$ is abundant
in transcendental numbers. But a transcendental number reported in
$\mathcal{C}_3$ is, to the knowledge of authors, only the
Komornik-Loreti constant $\delta=1.787231650\cdots$. See
\cite{KL}.

This section contains the proof that all Sturmian numbers are
transcendental. That enriches $\mathcal{C}_3$ with transcendental
numbers of continuum cardinality. In fact Sturmian words hitherto
have given birth to transcendental numbers in other manners.
Ferenczi and Mauduit showed real numbers whose expansions in some
integer base are Sturmian are transcendental \cite{FM}. Moreover,
they also generalized the transcendence results for Arnoux-Rauzy
sequences (on $3$ letters), and later it was extended to
Arnoux-Rauzy sequences on any $k$ letters by Risley and Zamboni
\cite{RZ}. Recently it is known that if the sequence of partial
quotients of the continued fraction expansion of a positive real
number is Sturmian, then the number is transcendental \cite{ADQZ}.
 For our purpose we need a classical result on transcendence of
functions in complex variables.

\begin{proposition}[\cite{Ma,LP}]\label{transcendence}
Let a function $f$ be defined by
$$f(w,z)=\sum_{n=1}^\infty \lfloor n w\rfloor z^n ,$$
where $w$ is real and $z$ is complex with $|z|<1$. Then
$f(\omega,\alpha)$ is transcendental if $\omega$ is irrational and
$\alpha$ is a nonzero algebraic number with $|\alpha|<1$.
\end{proposition}

Indeed Mahler \cite{Ma} proved in 1929 the preceding result for
quadratic irrational $\omega$, and later Loxton and van der
Poorten \cite{LP} extended the case to an arbitrary irrational
$\omega$. We are now in a position to state the main result of
this section.

\begin{theorem}
Every Sturmian number is transcendental, that is, if $d_\beta (1)$
is Sturmian, then $\beta$ is transcendental.
\end{theorem}

\proc{Proof.} Suppose $\mathrm{alph}(d_\beta (1))=\{a,b\}$ with
$0\leq a<b=\lfloor\beta\rfloor$. Since $\beta$ is Sturmian we have
for some irrational $\alpha \in(0,1)$,
\begin{eqnarray*}
1-\frac{b-a}{\beta}&=&\sum_{n=0}^\infty \frac{(b-a)S_{\alpha,0}
(n) + a}{\beta^{n+1}}
  =\sum_{n=0}^\infty \frac{(b-a)(\lfloor\alpha(n+1)\rfloor-\lfloor\alpha n\rfloor)+
  a}{\beta^{n+1}}\\
  &=&(b-a) \sum_{n=0}^\infty \frac{\lfloor\alpha(n+1)\rfloor-\lfloor\alpha n \rfloor}{\beta^{n+1}}
   +\frac{a}{\beta-1} .
\end{eqnarray*}
Thus the following equality holds.
\begin{equation}\label{contradiction}
\sum_{n=0}^\infty \frac{\lfloor\alpha(n+1)\rfloor-\lfloor\alpha n
\rfloor}{\beta^{n+1}}
  = \frac{1}{b-a}\left( 1-\frac{b-a}{\beta}-\frac{a}{\beta-1}
  \right).
\end{equation}
On the other hand,
$$\sum_{n=0}^\infty \frac{\lfloor\alpha(n+1)\rfloor-\lfloor\alpha n
\rfloor}{\beta^{n+1}}
  =\sum_{n=1}^\infty \frac{\lfloor\alpha n\rfloor}{\beta^n}
   -\frac{1}{\beta} \sum_{n=1}^\infty \frac{\lfloor\alpha
   n\rfloor}{\beta^n}
  =\left(1-\frac{1}{\beta}\right) \sum_{n=1}^\infty \frac{\lfloor\alpha
  n\rfloor}{\beta^n}.$$
If $\beta$ were algebraic, the left side of Equation
(\ref{contradiction}) would be transcendental by Proposition
\ref{transcendence} whereas the right one algebraic. \ep \medbreak

\proc{Example 5.} We assume
$d_\beta(1)=313113131131131311313113113\cdots$, where the sequence
is obtained by substituting $3$ for $1$ and $1$ for $0$ in the
sequence $1c_{\tau^{-2}}$ of Example 3. Then such $\beta\in(3,4)$
exists and is transcendental. Furthermore, $\mathrm{diam}
\{T_\beta ^n 1\}_{n\geq0}=2/\beta$ and $\overline{\{T_\beta^n
1\}}_{n \geq 0}$ is of Lebesgue measure zero. \medbreak

%
%
%
%       R  E  F  E  R  E  N  C  E  S
%
%
%

\hfill\\

\noindent School of Mathematical Sciences,\\
Seoul National University,\\
Seoul 151-747, Korea.\\
E-mail: dpchi@math.snu.ac.kr\\ \vspace{0.3cm}

\noindent School of Computational Sciences,\\
Korea Institute for Advanced Study,\\
Seoul 130-722, Korea.\\
E-mail: doyong@kias.re.kr

\begin{thebibliography}{30}
\bibitem{AC1} J.-P. Allouche and M. Cosnard.
The Komornik-Loreti constant is transcendental. \textit{Amer.
Math. Monthly} \textbf{107} (2000) 448-449.

\bibitem{AC2} J.-P. Allouche and M. Cosnard.
Non-integer bases, iteration of continuous real maps, and an
arithmetic self-similar set. \textit{Acta Math. Hungar.}
\textbf{91} (2001) 325-332.

\bibitem{ADQZ} J.-P. Allouche, J.L. Davison, M. Queff\'{e}lec and L.Q. Zamboni.
Transcendence of Sturmian or morphic continued fractions.
\textit{J. Number Theory} \textbf{91} (2001) 39-66.

\bibitem{BM} A. Bertrand-Mathis.
D\'{e}veloppement en base $\theta$ et r\'{e}partition modulo 1 de
la suite $(x\theta^n )$. \textit{Bull. Soc. Math. France}
\textbf{114} (1986) 271-324.

\bibitem{Bl} F. Blanchard.
$\beta$-expansions and symbolic dynamics. \textit{Theoret. Comput.
Sci.} \textbf{65} (1989) 131-141.

\bibitem{BL} J.-P. Borel and F. Laubie.
Quelques mots sur la droite projective r\`{e}elle. \textit{J.
Th\'{e}or. Nombres Bordeaux} \textbf{5} (1993) 23-51.

\bibitem{CH} E.M. Coven and G.A. Hedlund.
Sequences with minimal block growth. \textit{Math. Systems Theory}
\textbf{7} (1973) 138-153.

\bibitem{FM} S. Ferenczi and C. Mauduit.
Transcendence of numbers with a low complexity expansion.
\textit{J. Number Theory} \textbf{67} (1997) 146-161.

\bibitem{He} G.A. Hedlund.
Sturmian minimal sets. \textit{Amer. J. Math.} \textbf{66} (1944)
605-620.

\bibitem{Ho} F. Hofbauer.
$\beta$-shifts have unique maximal measure. \textit{Monatsh.
Math.} \textbf{85} (1978) 189-198.

\bibitem{KL} V. Komornik and P. Loreti.
Unique developments in non-integer bases. \textit{Amer. Math.
Monthly} \textbf{105} (1998) 636-639.

\bibitem{Lo} M. Lothaire.
\textit{Algebraic combinatorics on words.} Cambridge University
Press, 2002.

\bibitem{LP} J.H. Loxton and A.J. van der Poorten.
Arithmetic properties of certain functions in several variables
III. \textit{Bull. Austral. Math. Soc.} \textbf{16} (1977) 15-47.

\bibitem{Ma} K. Mahler.
Arithmetische Eigenschaften der L\"{o}sungen einer Klasse von
Funktionalgleichungen. \textit{Math. Ann.} \textbf{101} (1929)
342-366.

\bibitem{Mi} F. Mignosi.
Infinite words with linear subword complexity. \textit{Theoret.
Comput. Sci.} \textbf{65} (1989) 221-242.

\bibitem{MH} M. Morse and G.A. Hedlund.
Symbolic dynamics II: Sturmian sequences. \textit{Amer. J. Math.}
\textbf{62} (1940) 1-42.

\bibitem{Pa} W. Parry.
On the $\beta$-expansion of real numbers. \textit{Acta Math. Acad.
Sci. Hungar.} \textbf{11} (1960) 401-416.

\bibitem{Pv} B. Parvaix.
Propri\'{e}t\'{e}s d'invariance des mots sturmiens. \textit{J.
Th\'{e}or. Nombres Bordeaux} \textbf{9} (1997) 351-369.

\bibitem{RZ} R.N. Risley and L.Q. Zamboni.
A generalization of Sturmian sequences: combinatorial structure
and transcendence. \textit{Acta Arith.} \textbf{95} (2000)
167-184.

\bibitem{Sc} J. Schmeling.
Symbolic dynamics for $\beta$-shift and self-normal numbers.
\textit{Ergodic Theory Dynam. Systems} \textbf{17} (1997) 675-694.


\end{thebibliography}
\end{document}